\documentstyle[12pt]{article}
\newcommand{\R}{\mbox{I$\!$R}}

\newcommand{\qed}{{\hfill {$\rlap{$\sqcap$}\sqcup$}}\\[0.2in]\hspace*{0.5in}}
\newcommand{\qedwh}{{\hfill {$\rlap{$\sqcap$}\sqcup$}}\\[0.2in]}
\newcommand{\bk}{\\[0.03in] \hspace*{0.5in} }
\newcommand{\btd}{\bigtriangledown}

\newcommand{\mfor}{\ \ \ \ {\mbox{for}} \ \ }
\newcommand{\for}{\ \ \ \ for \ \ }

\begin{document}

\vspace*{0.07in}
\begin{center} 
{\LARGE   {\bf Asymptotic Behavior of }}
\medskip \\ {\LARGE {\bf  Positive Solutions of the Equation }}
\medskip \\ {\LARGE {\bf $\Delta u + K u^{{n + 2}\over {n - 2}} = 0$  
in ${\R}^n$}}
\medskip \medskip\\ {\LARGE {\bf and Positive Scalar Curvature}}
\medskip \medskip \medskip  \medskip\\ 
{\large { Ka-Luen {\Large C}HEUNG\footnote{Department of Mathematics, 
Rutgers University, Hill Center-Busch Campus,\\
110 Frelinghuygen Rd.,
 Piscataway, NJ 08854, U.S.A. \ \ \ \ {\tt cheung@math.rutgers.edu}}
 \ \ and \ \  Man-Chun {\Large L}EUNG\footnote{
Department of
Mathematics,  National University of Singapore, Lower Kent Ridge Rd.,\\
Singapore 119260, Republic of Singapore \ \
{\tt matlmc@math.nus.edu.sg} }}}
\end{center}

\vspace*{0.15in} 

\begin{abstract}We study asymptotic behavior of positive smooth
solutions of the conformal scalar curvature equation in ${\R}^n$. We
consider the case when the scalar curvature of the conformal metric is
bounded between two positive numbers outside a compact set.  
It is shown that the solution  has slow decay if the radial change is 
controlled.  For a positive solution with slow decay, the corresponding conformal metric is found to 
be complete if and only if the total volume is infinite. 
We also determine the sign of the Pohozaev
number in some situations and  show that if the Pohozaev is equal to zero, then either the solution has fast decay, or the conformal metric corresponding to the 
solution is complete and the corresponding solution in $\R \times S^{n - 1}$ 
has a sequence of local maxima 
that approach the standard
spherical solution. 
   
\end{abstract}

\vspace*{0.1in} 

1991 AMS MS Classifications: Primary 35J60\,; Secondary 58G03.

\vspace*{0.5in}

{\bf \Large {\bf 1. \ \ Introduction}}

\vspace*{0.3in}

For $n \ge 3$, we consider positive smooth solutions $u$ of the 
conformal scalar curvature equation
$$
\Delta u + K u^{{n + 2}\over {n - 2}} = 0 \ \ \ \ {\mbox{in}} \ \ {\R}^n, \leqno (1.1)
$$
where $\Delta$ is the standard Laplacian  on ${\R}^n$ and $K$
is a smooth function.  Equation (1.1) is studied extensively
by many authors and  remarkable results are obtained. Some of the 
 recent publications \cite{Bianchi} \cite{Caffarelli-Gidas-Spruck}
\cite{Chen-Lin.2} \cite{Chen-Lin.3} \cite{Egnell}
\cite{Gidas-Ni-Nirenberg.2} \cite{Gidas-Spruck} \cite{K-M-P-S}
\cite{Y.Y.Li.1} \cite{Y.Y.Li.2} \cite{Y.Li} \cite{Li-Ni-2}
\cite{Li-Ni-3} \cite{Lin.1} \cite{Mazzeo-Pacard} 
\cite{Mazzeo-Pollack-Uhlenbeck} \cite{Schoen}
\cite{Zou.1} 
 on this equation provide excellent  
descriptions of the development.  
Here we are 
 concerned with the case when $K$ is bounded
between two positive numbers outside a compact subset of
${\R}^n$. That means that the conformal metric $g = u^{4/(n - 2)} g_o$
has bounded positive scalar curvature outside a compact subset, where
$g_o$ is Euclidean metric. Accordingly, positive smooth solutions of
equation (1.1) can be broadly classified as whether the conformal
metric $g$ is complete or not, and if it is not complete, whether it can
be realized as a smooth metric on $S^n$ after a stereographic
projection.\bk 
The total volume of the Riemannian manifold 
$({\R}^n,\, g)$ is given by
$$
\int_{{\R}^n} u^{{2n}\over {n - 2}} (x) \, dx\,. \leqno (1.2)
$$
One can also classify  positive smooth solutions of
equation (1.1) according to whether the total volume is finite or not. These geometric
classifications are closely related to the decay of $u$ near infinity. A  positive smooth solution
$u$ of equation (1.1) in ${\R}^n$ is said to have {\it slow decay}
\cite{Zou.2} 
if there
exist positive constants $C_s$ and $r_s$ such that 
$$
u (x) \le C_s |x|^{ - {{n - 2}\over 2} } \mfor |x| \ge r_s\,. \leqno (1.3)
$$
It is said to have {\it fast decay} if 
there exist positive constants $C_f$ and $r_f$ such that 
$$
u (x) \le C_f |x|^{ - (n - 2) } \mfor |x| \ge r_f\,. \leqno (1.4)
$$
(Note that the definition of slow decay is different from the one in 
\cite{Chen-Lin.1} and \cite{Chen-Lin.2}.) If $u$  has fast decay, then $g$
can be realized as a metric on $S^n$ (cf. lemma 2 in \cite{K-M-P-S}, see also
\cite{Chen-Lin.1}). 
Using a result
of Brezis and Kato \cite{Brezis-Kato}, it is shown in \cite{Leung.5} 
that if the
total volume of $({\R}^n, g)$ is finite, then $u$ has fast decay.\bk
We may infer 
that the cases of slow decay, completeness and infinite total volume
are related to each other. They also 
pertain to  the asymptotic dimension of $({\R}^n, \,g)$ (cf. the
conjecture propounded by Gromov in \cite{Gromov}, see also \cite{Leung.5}). For slow
decay means that the conformal metric $g$ is asymptotically a
bounded perturbation of the standard product metric on ${\R}^+ \times
S^{n - 1}$ \cite{Leung.5}. When $K$ is equal to a positive constant outside a compact
set, Caffarelli, Gidas and Spruck in their celebrated work
\cite{Caffarelli-Gidas-Spruck} 
show that either $g$  can be            
realized as a smooth metric on $S^n$, or $u$ is asymptotic to one
of a one-parameter family of radial solutions and has slow decay.  
In \cite{Chen-Lin.1}, using the method of moving planes,  Chen and Lin show that if there exist positive constants
$c_1$, $c_2$, $c_3,$ $l$ and $R_o\,,$ with $c_3 > 1$, such that 
$$
0 < c_1 \le | \btd K (x) | \, |x|^{l + 1} \le c_2 \mfor |x| \ge R_o\,,
\leqno (1.5)
$$
$$
\lim_{|x| \to \infty} K (x) = K_\infty > 0  \ \ \ \ {\mbox{and}} \leqno (1.6)
$$ 
$$
K (y) \le K (x) \mfor |y| \ge c_3 \,|x| \ \ \ \ {\mbox{and}} \ \ |x| \ge
R_o\,, \leqno (1.7)
$$
then $u$ has slow decay.  In
\cite{Leung.5}, 
the second author
shows that if $K$ is bounded between two positive numbers in ${\R}^n$, $x \cdot \btd K (x) \le 0$ for 
large $|x|$, $u$ is bounded from above 
and $u$ is {\it radially dominating} (see \cite{Leung.5}), then $u$ has slow
decay. In this paper we obtain the slow decay by assuming a control on
the growth in the radial direction.\\[0.2in]
{\bf Theorem A.} \ \ {\it Assume that $K$ is bounded between two
positive constants outside a compact subset of ${\R}^n$. Assume also that there
exist positive constants $r_o$, $c$ and  $C$  such that}
$$
{{\partial K}\over {\partial r}} (r, \theta) \ge - C e^{-c\,r} \for
\theta \in S^{n - 1} \ \ \ \ and \ \ r \ge r_o\,,\leqno (1.8)
$$
{\it where $(r, \theta)$ is the polar coordinates on ${\R}^n$. Let $u$ be a positive smooth solution of equation (1.1) in
${\R}^n$. If there
exist positive constants $C'$ and $r'$ such that}
$$
r  \bigg\vert{{\partial u}\over {\partial r}} 
(r, \theta)
\bigg\vert \le C'\, u (r, \theta) \for \theta \in S^{n - 1} \ \ \ \
and \ \ \ \ r \ge r'\,, \leqno (1.9)
$$
{\it then $u$ has slow decay.}\\[0.2in]
\hspace*{0.5in}It is observed that if there exist positive 
constants $C_l$ and $r_l$ such that 
$$
u (x) \ge C_l\, |x|^{ - {{n - 2}\over 2} } \mfor |x| \ge r_l\,, \leqno (1.10)
$$
then $g$ is a complete metric on ${\R}^n$. Inequality (1.10) is related
to the {\it Pohozaev number} of $u$, denoted by $P (u)$ (see section
2). For a positive smooth solution $u$ with slow decay, $P (u) \not= 0$ implies that (1.10)
holds. Using the results of Korevaar, Mazzeo, Pacard and Schoen in 
\cite{K-M-P-S}, if 
$K$ is equal to a positive constant outside a compact set and $g$ is 
complete, then (1.10) holds and $P (u) < 0$. Chen and Lin
\cite{Chen-Lin.1} 
show that,
under the assumptions (1.5), (1.6), (1.7) and in addition, if it is also assumed that $x \cdot \btd K (x) \le 0$ for large
$|x|$,  then the conformal metric $g$ can be realized
as a smooth metric on $S^n$ if and only if $P (u) = 0$, and it 
can be deduced  that  $g$ is
complete (or (1.10) holds) if and only if $P (u) < 0$. 
The sign of the Pohozaev number can also be determined in the following situation.
\\[0.2in]
%\pagebreak
{\bf Theorem B.} \ \ {\it  Assume that $\lim_{|x| \to \infty} K (x) =
K_\infty > 0$ and $|\btd
K| \le C_o$ in ${\R}^n$ for a positive number $C_o$. Let $u$ be
a positive smooth solution of (1.1) in ${\R}^n$ with slow decay. Assume
also that $K$ satisfies one of the conditions in lemma 2.4 so that $P
(u)$ exists. Then we have $P (u) \le 0$. In addition, if it is also assumed that $x \cdot
\btd K (x) \le 0$ for large $|x|$, then $P (u) = 0$ implies that $u$ has fast decay.}\\[0.2in]
\hspace*{0.5in}Similar results are obtained by Chen and Lin in \cite{Chen-Lin.1} and
\cite{Chen-Lin.4}, and by Korevaar, Mazzeo, Pacard and Schoen in \cite{K-M-P-S}. Indeed, our
contribution to theorem B consists mainly of showing that the arguments used in their works can be
applied in a general setting. It is rather natural to relate the condition that 
$P (u) = 0$ with fast decay, due to the results mentioned above, and also to the Kazdan-Warner
identity which says that if the conformal metric
$u^{4/(n - 2)} g_o$ can be realized as a smooth metric on $S^n$,  then $P (u) = 0$. An interesting 
example of Chen and Lin (theorem 1.6 in 
\cite{Chen-Lin.4}) shows that the relation is actually quite subtle. Let 
$$
v (s, \theta) = r^{{n - 2}\over 2} u (r, \theta)\,, \ \ \ \ r = e^s\,,
\ \ r > 0\,, \ \  s \in \R  \ \ {\mbox{and}} \ \ \theta \in S^{n -
1}. \leqno (1.11)
$$
We 
show that, under the conditions in the first part of theorem B, if $P (u) = 0$, then either $u$ has
fast decay, or   the conformal metric $g$ is complete and 
$$
\liminf_{|x| \to \infty} |x|^{{n - 2}\over 2} u (x) = 0
$$
and  
there exists a sequence $\{s_j\} \subset \R$ of local maxima of 
$$
\bar v = \int_{S^{n - 1}} v \, d\theta\,, \leqno (1.12)
$$
$s_j \to +\infty$ as $j \to \infty\,,$ such that the 
sequence $\{v_j\}$ defined by
$$
v_j (s, \,\theta) = v (s_j + s, \,\theta) \mfor s \in \R\,, \ \ \theta \in S^{n - 1} \ \ {\mbox{and}} \ \ 
j = 1, \ 2,..,
$$
converges uniformly in $C^2$-norm on compact subsets of $\R \times S^{n - 1}$ 
to $C \,(\cosh s)^{(2 - n)/2}$ for a positive constant $C$ depending on $K_\infty$.  
The example in \cite{Chen-Lin.4} suggests that the second situation in the above statement may occur. 
We note that the conformal factor $(\cosh s)^{-2}$ transforms the cylinder 
$\R \times S^{n - 1}$ into $S^n \setminus \{p, \,- p\}$ \cite{K-M-P-S}. It is also noteworthy that even if
$P (u) = 0$, the  conformal metric $g$ may still be complete.\\[0.2in]
{\bf Theorem C.} \ \ {\it Assume that $K$ is
 bounded between two
positive numbers outside a compact subset of ${\R}^n$. Let $u$ be
a positive smooth solution of equation (1.1) in ${\R}^n$ with slow decay. Then either $u$ has fast 
decay or the conformal metric $g = u^{4/(n - 2)} g_o$ is complete. Furthermore,  
the conformal metric $g$ is complete  if and only if the total volume of $({\R}^n, g)$ is
infinite.}\\[0.2in]
\hspace*{0.5in}When $K$ may not have a 
limit at  infinity, for a particular class of positive
solutions, we can determine the sign of the Pohozaev number.
%\\[0.2in]

\pagebreak

{\bf Theorem D.} \ \ {\it Assume that $K$ is
 bounded between two
positive numbers outside a compact subset of ${\R}^n$. Let $u$ be a
positive smooth solution of equation (1.1) in ${\R}^n$. Assume that }
$$
\lim_{r \to \infty}  \int_{S_r} r \left[ {{\partial u}\over {\partial r}} + {{n
- 2}\over 2} {u\over r} \right]^2 \, dS = 0\,. \leqno (1.13)
$$
{\it If $P (u)$ exists, then $P (u) \le 0$. Furthermore, $P (u) = 0$
implies that}
$$
\liminf_{|x| \to \infty} |x|^{(n - 2)/2}\, u (x) = 0\,.
$$

\vspace*{0.1in}
                                           
\hspace*{0.5in}In
theorem D, we do not assume that $u$ has slow decay. The meaning of condition (1.13) is explained in
more detail in  (2.26) and section 4. We outline  the content of each section. Section 2
reviews the Pohozaev identity, gradient estimate and spherical Harnack
inequality. Theorem A is proved in section 3. Section 4 is
devoted to a study of the Pohozaev number in general and 
the proofs of theorem B, C and D in particular. We use $C$,
$C_1,...,$ $c$, $c_1,...$ to denote various constants which may be
different from section to section according to the context. Moreover, throughout this article, 
$n \ge 3$ is an integer.

\vspace{0.5in}

%\pagebreak         

{\bf \Large {\bf 2. \ \ Preliminaries}}

\vspace{0.3in}

Let $u$ be a positive smooth solution of equation (1.1) and 
$$
P (u, r) = {{n - 2}\over {2n}} \int_{B_o (r)} x \cdot \btd K (x) \, u^{{2n}\over {n - 2}}
(x) \, dx \mfor r > 0\,, \leqno (2.1)
$$
where $B_o (r)$ is the open ball with center at the origin and radius 
 $r > 0$. 
The Pohozaev identity \cite{Pohozaev} \cite{Ding-Ni} shows that 
$$
P (u, r) =  \int_{S_r} \left[ r \left( {{\partial u}\over {\partial r}} \right)^2
 - {r\over 2} |\btd u|^2 +  {{n - 2}\over {2n}} r  K u^{{2n}\over {n -
2}}  + {{n - 2}\over 2} \,u \,{{\partial
u}\over {\partial r}} \right] \,dS  \leqno (2.2)
$$
for $r > 0$, where $S_r = \partial B_o (r)$ is the sphere of radius
$r$. The Pohozaev number of $u$ is given by 
$$
P (u) = \lim_{r \to \infty} P (u, r) \leqno (2.3)
$$
provided the limit exists.
%\\[0.2in]

\pagebreak

{\bf Lemma 2.4.} \ \ {\it  Let $u$ be a positive smooth solution of
equation 
(1.1)
with 
slow decay. Then the limit in (2.3) exists  if we assume 
anyone of the following conditions.}\\[0.05in]
({\bf I}) \ \ {\it There exists a  number $m > 1$ such that} 
$$
x \cdot \btd K = r{{\partial K}\over {\partial r}} \in L^m \left({\R}^n \setminus 
B_o (1) \right)\,. \leqno (2.5)
$$
({\bf II}) \ \ {\it There exist positive numbers $C$, $\varepsilon$ and $r_o$ such that}
$$
\bigg\vert {{\partial K}\over {\partial r}} \bigg\vert \le {C\over {r
(\ln \, r)^{1 + \varepsilon} }} \for  r \ge r_o\,. \leqno (2.6)
$$
({\bf III}) \ \ {\it $x \cdot \btd K (x)$ does not change sign for large
$|x|$ and there exist
positive constants $\alpha$ and $\beta$ such that $- \alpha^2 \le K
\le \beta^2$ in ${\R}^n$.}\\[0.1in]
{\bf Proof.} \ \ We note that 
$$
\bigg\vert \int_{ {\R}^n \setminus B_o (R) } |x| {{\partial K}\over
{\partial r}} u^{{2n}\over {n - 2}} \, dx \bigg\vert \le C_o 
\left[  \int_{ {\R}^n \setminus B_o (R) } \bigg\vert r{{\partial K}\over
{\partial r}} \bigg\vert^m \, dx \right]^{1\over m} 
\left[  \int_{ {\R}^n \setminus B_o (R) } |x|^{- nl} 
\, dx \right]^{1\over l} \to 0
$$
as $R \to \infty$, where $l > 1$ is the number such that $1/m + 1/l =
1$. Here $C_o$ is a positive constant. Likewise,  
$$
\bigg\vert \int_{ {\R}^n \setminus B_o (R) } |x| {{\partial K}\over
{\partial r}} u^{{2n}\over {n - 2}} \, dx \bigg\vert \le C'
\int_R^\infty {{dr}\over {r (\ln \, r)^{1 + \varepsilon} }} \to 0
$$
as $R \to \infty$, where $C'$ is a positive constant. Hence the limit in (2.3) exists if we assume either
(I) or (II). Assume that $x \cdot \btd K (x) \le 0$ for large $|x|$.  
The integral in (2.1) 
is a decreasing function for large $r$. For $R > R_o$ large enough, we
have
\begin{eqnarray*} 
\int_{R_o}^R \int_{S_r} r {{\partial K}\over {\partial r}}
u^{{2n}\over {n - 2}} \, dS \, dr 
& \ge &  C_1 \int_{R_o}^R \int_{S^{n-1}} {{\partial K}\over {\partial r}}
\,d \theta \, dr
 =  C_1 \int_{R_o}^R  {d\over {dr}} \left( \int_{S^{n - 1}} K (r, \theta)\,
d\theta \right) \, dr\\ & \ge & -C_1\, \omega_n\, (\alpha^2 + \beta^2)\,, 
\end{eqnarray*}
where $C_1$ is a positive constant and $\omega_n$ is the volume of
$S^{n - 1}$. Hence the integral in (2.1) is bounded from below and the limit in
(2.3) exists. The other case is similar.\qed
We assume that there exist positive constants $a$, $b$ and $r_o$ such
that 
$$
0 < a^2 \le K (x) \le b^2 \mfor |x| \ge r_o\,. \leqno (2.7)
$$
$u$ is said to satisfy a spherical Harnack inequality if there exists
a positive constant $C_h$  such that the inequality 
$$
\sup_{S_r} u \, \le  \,C_h\, \inf_{S_r} u  \leqno (2.8)
$$
holds for all $r > 0$. 
In \cite{Lin.1}, the following relations are observed (see also \cite{Zou.2}).  A version of the
result is proved in
\cite{Chen-Lin.4}. It is certainly well-known among the experts, but for the sake of
completeness, we describe the argument below.\\[0.2in] {\bf Lemma 2.9.} \ \ {\it Assume
that
$K$ satisfies (2.7). 
 Let $u$ be a positive smooth solution of (1.1)
in ${\R}^n$. Then the following statements are equivalent.\\[0.05in]
(a) \ \ $u$ has slow decay.\\[0.05in]
(b) \ \ $u$ satisfies a spherical Harnack inequality.\\[0.05in]
(c) \ \ There exists a positive constant $C_g$  such that $r
|\btd u (r, \theta)| \le C_g\, u (r, \theta)$ for $r > 0$ and $\theta \in
S^{n - 1}$.}\\[0.1in]
{\bf Proof.} \ \ Assuming slow decay, 
write equation (1.1) as $\,\Delta u (x) + f (x)\, u
(x) = 0$ in ${\R}^n$, where $f (x) = K (x)\, u^{4/(n - 2)} (x)$. 
We have the following scaling property for $f$:
$$
0 \le R^2 f^2 (R x) \le {{b^2\, C_s^{4\over {n - 2}} }\over {|x|^2}} \mfor
R\, |x| \ge r_1\,,
$$
where $r_1 = \max \ \{ r_s, r_o\}$. Here $C_s\,, r_s$ are the constants 
in (1.3) and $r_o$ is the constant in (2.7). 
Using a scaling as above and the 
Harnack inequality (theorem 8.20 in \cite{Gilbarg-Trudinger}), there
exist positive constants $C_1$ and $R_1$  such that 
$$
\sup_{B_{2 R}\setminus B_R} u  \le C_1 \,  \inf_{B_{2 R}\setminus B_R} u
\mfor R \ge R_1\,. \leqno (2.10)
$$
It follows that we have a spherical Harnack inequality for spheres of
large radius. As $u$ is positive, therefore inside a big ball $B_o
(R_2)$ we have 
$$\sup_{B_o (R_2)} u \le C_2 \inf_{B_o (R_2)} u \leqno (2.11)$$
for a positive constant $C_2$ that depends on $R_2$ and $u$. 
Hence we have
(b). We obtain (c)
by using  the gradient estimate in \cite{Gilbarg-Trudinger} (pp. 37), 
equation (1.1) and (2.8). On the other hand, (c) implies that, for a
fixed $r > 0$,  
$$
|\btd_\theta u (r, \theta) | \le C' \,u (r, \theta) \mfor \theta \in S^{n
- 1}\,,
$$
where $C'$ is a positive constant independent on $r$. After an integration we have 
$$
\sup_{S_r} u  \le e^{C' \pi}\, \inf_{S_r} u  \mfor r > 0\,,
$$
which is (b). From (b) and the fact that there exist positive
constants $C_2$ and $r_2$ such that 
$$
\int_{S^{n - 1}} u (r, \theta) \, d\theta \le C_2 \,r^{- (n - 2)/2}
\mfor r \ge r_2 \leqno (2.12)
$$
(see, for example, \cite{Leung.5}), we conclude that $u$ has slow decay.\qed
The following result is a direct consequence of Pohozaev identity (2.2) and
lemma 2.9.\\[0.2in]
{\bf Lemma 2.13.} \ \ {\it Assume that $K$ satisfies (2.7).   
Let $u$ be a positive smooth solution of
equation (1.1)
in ${\R}^n$ with slow decay. Assume also that $K$ satisfies one of the
conditions in lemma 2.4 so that $P (u)$ exists. If
$P (u) \not= 0$,  then 
$$u (x) \ge C_l |x|^{-(n-2)/2} \leqno (2.14)$$
 for large $|x|$ and for a
positive constant $C_l$.}\\[0.1in]
{\bf Proof.} \ \ Suppose that (2.14) does not hold. Then there exists
a sequence  $\{ x_i \} \subset {\R}^n$ 
such that $\lim_{i \to \infty} |x_i| = \infty$ and $u (x_i) |x_i|^{(n
- 2)/2} \to 0$ as $i \to \infty$. Let $r_i = |x_i|$ for $i = 1, 2,...$ 
Using spherical Harnack inequality (2.8) and the gradient estimate in
lemma 2.9  we have 
$$
r_i^{(n - 2)/2} \,[ \,\max_{S_{r_i}} \, u\,] \to 0 \ \ \ \ {\mbox{and}} \
\ \ \  
 r_i^{n/2} \, [\, \max_{S_{r_i}} \, | \btd  u| \,] \to 0 \ \ \ \
{\mbox{as}} \ \ i \to \infty\,.
$$ 
From (2.2) we have $P (u) = 0$, which is a contradiction. Hence (2.14) holds.\qed
Without the assumption of slow decay, we have the
following estimates.\\[0.2in]
{\bf Theorem 2.15.} \ \ {\it Assume that $K$ satisfies (2.7). Let $u$ be
a positive smooth solution of equation (1.1)
in ${\R}^n$. If there exists a positive constant $\delta$ such that}
$$
| P (u, r) | \ge \delta^2 \leqno (2.16)
$$
{\it for large $r$, then there exist positive constants $C'$ and $C''$
such that}
$$
\int_{B_o (r)} u^{{2n}\over {n - 2}} (x) \,dx \ge C' \ln r \ \ \ \ and \
\ \ \ 
\int_{B_o (r)} |\btd u|^2 \,dx \ge C'' \ln r \leqno (2.17)
$$
{\it for large $r$. In particular, if $P (u)$ exists and is non-zero, 
then the total volume of $({\R}^n, \,g)$ is infinite.}\\[0.1in]
{\bf Proof.} \ \  Applying Young's inequality we have
\begin{eqnarray*}
(2.18) \ \ \ \ \  \ \ \ \ \ \int_{S_r} u \,{{\partial
u}\over {\partial r}}  \,dS & \le & \int_{S_r} r \bigg\vert {{\partial
u}\over {\partial r}} \bigg\vert^2  \,dS  +  \int_{S_r} {{u^2}\over r} \, dS\\
& \le & \int_{S_r} r | \btd u |^2 \, dS + C_1 \int_{S_r} r \,u^{{2n}\over {n -
2}} \,dS + {{\delta^2}\over {n - 2}}  \ \ \ \ \  \ \ \ \ \  \ \ \ \ \ 
\end{eqnarray*}
for large $r$, where $C_1$ is a positive constant depending on $\delta$ and $n$. From 
Pohozaev identity (2.2) and (2.18), there exist positive constant $C_2$ and $C_3$ such that
$$
C_2 \int_{S_r} | \btd u |^2 \, dS  + 
C_3 \int_{S_r} u^{{2n}\over {n -
2}} \,dS  \,\ge \,{{\delta^2}\over {2r}} \leqno (2.19)
$$
for large $r$. It follows that 
$$
\int_{B_o (r)} | \btd u |^2 \, dx + \int_{B_o (r)} u^{{2n}\over {n -
2}} \, dx \ge C_4 \ln r \leqno (2.20)
$$
for large $r$ and for a positive constant $C_4$. We can modify the argument in the proof of
theorem 3.1 in \cite{Leung.5} for our case so as to obtain positive constants 
$C_5$ and $D_1$ such that 
$$
C_5 \int_{B_o (2R)} u^{{2n}\over {n -
2}} \, dx + D_1 \ge \int_{B_o (R)} | \btd u |^2 \, dx \leqno (2.21)
$$
for all large $R$. (2.20) and (2.21) imply that 
$$
\int_{B_o (r)} u^{{2n}\over {n -
2}} \, dx \ge C' \ln r \leqno (2.22)
$$
for large $r$ and for a positive constant $C'$. Similarly we obtain the other estimate in
(2.17).\qedwh 
{\bf Remark 2.23.} \ \ Partial results regarding the estimates in (2.17) are obtained
in 
\cite{Ding-Ni} and \cite{Leung.5}.  
The arguments in \cite{Ding-Ni} and \cite{Leung.5} yield the stronger estimates that
$$
\int_{S_r} u^{{2n}\over {n - 2}} \,dS \ge {{C_6}\over {r}} \ \ \ \ {\mbox{and}} \ \ \ \ 
\int_{S_r} |\btd u|^2 \,dS \ge {{C_7}\over {r}}
$$
for large $r$ and for some positive constants $C_6$ and $C_7$.\\[0.2in]
\hspace*{0.5in}Let $v$ be defined as in (1.11). 
Then $v$ satisfies the equation 
$$
{{\partial^2 v}\over {\partial s^2}} + \Delta_\theta\, v - \left( {{n -
2}\over 2} \right)^2 v + K v^{{n + 2}\over {n - 2}} = 0 \ \ \ \
{\mbox{in}} \ \ \R \times S^{n - 1}, \leqno (2.24)
$$ 
where $\Delta_\theta$ is the standard Laplacian on $S^{n - 1}$ (see,
for example, \cite{Leung.5}). The
Pohozaev identity (2.2) becomes
\begin{eqnarray*}
(2.25) \  P (v, s) &  = & {{n - 2}\over {2n}} \int_{-\infty}^s \int_{S^{n - 1}} {{\partial K}\over
{\partial t}} v^{{2n}\over {n - 2}} \, d\theta \, dt\\
&  =  & \int_{S^{n - 1}} \left[ {1\over 2} 
\left( {{\partial v}\over {\partial s}} \right)^2  
-  {1\over 2} | \btd_\theta v |^2 
- {1\over 2} \left( {{n - 2}\over 2} \right)^2  v^2
 + {{n - 2}\over {2n}}  K  v^{{2n}\over {n
- 2}} \right] \, d\theta \ \ \ \ \ 
\end{eqnarray*}
for $s \in \R$. With the notations, condition (1.13) in theorem D can
also be written as
$$
\lim_{s \to \infty} \int_{S^{n - 1}} \left[ {{\partial v}\over
{\partial s}} (s, \theta) \right]^2  \, d\theta = 0\,. \leqno (2.26)
$$

\pagebreak

%\vspace{0.5in}

{\bf \Large {\bf 3. \ \ Slow decay}}

\vspace{0.3in}

{\bf Proof of Theorem A.} \ \ 
By differentiating the function
$$
\int_{S_r} u^{{2n}\over {n - 2}} \, dS = \int_{S^{n - 1}} u^{{2n}\over {n - 2}} (r, \theta)\,
r^{n - 1} d\theta \mfor r > r_o
$$
and using (1.9), there exist positive constants $C_o$ and $m$ such that 
$$
\int_{S_r} u^{{2n}\over {n - 2}} \, dS \le C_o\, r^m \mfor r \ge
r_o\,. \leqno (3.1)
$$
By (1.8) and (3.1) we have 
$$
{{n - 2}\over {2n}}
\int_{B_o (R)}  
r {{\partial K}\over {\partial r}} u^{{2n}\over {n - 2}} \,dx \ge - C_1 - C_2 \int_{r_1}^R
e^{-cr} r^{m + 1} \,dr \ge -C_3 \mfor R \ge r_1\,,
$$
where $C_1$, $C_2$ and $C_3$ are positive constants and $r_1 \ge \max
\, \{ r_o\,, r' \}$. Together with
(1.9) and  Pohozaev identity (2.2) we obtain
$$
\int_{S_R} R |\btd u|^2 \, dS \le 2 \,C_3 + C_4 \int_{S_R} {{u^2}\over R} \,
dS +  {{n - 2}\over n} \int_{S_R} R K u^{{2n}\over {n - 2}} \, dS \leqno (3.2)
$$
for large $R$, where $C_4$ is a positive constant. By choosing
a larger $C_3$ if necessary, we may assume that (3.2) holds for all $R >
0$. Similarly, as $u > 0$ in ${\R}^n$, we may assume that (1.9) holds
for $r > 0$. 
Integrating both sides of (3.2) we obtain  
$$
\int_{B_o (R)} r |\btd u|^2 \, dx \le 2 C_3\, R + C_4 \int_{B_o (R)}
{{u^2}\over r} dx  \, +  {{n - 2}\over n} \int_{B_o (R)} r K u^{{2n}\over {n - 2}}
\, dx \leqno (3.3) 
$$
for $R > 0$.  Using equation (1.1) and (1.9) we have 
\begin{eqnarray*}
(3.4) \ \ \ \ & \ & \int_{B_o (R)} r K u^{{2n}\over {n - 2}} \, dx = \int_{B_o (R)}
r u \left(K u^{{n + 2}\over {n - 2}}
\right) \, dx = \int_{B_o (R)} (ru) (- \Delta u) \, dx \ \ \ \ \ \ \ \\
& = & \int_{B_o (R)} r | \btd u |^2 \, dx + \int_{B_o (R)} u {{\partial
u}\over {\partial r}} \, dx - \int_{S_R} R u {{\partial u}\over {\partial r}}
\, dS\\
& \le & \int_{B_o (R)} r | \btd u |^2 \, dx + C_5 \int_{B_o (R)}
{{u^2}\over r} \, dx + C_6\int_{S_R} u^2 \, dS
\end{eqnarray*}
for large $R$, where $C_5$ and $C_6$ are positive constants. 
It follows from (3.3) and (3.4) that 
\begin{eqnarray*}
(3.5) \ \ \ \ \ & \ & {2\over n}  \int_{B_o (R)} r K u^{{2n}\over {n -
2}} \, 
dx\\ & \le & 2 C_3 R  + 
C_7\int_{B_o (R)} {{u^2}\over r} \, dx + C_8 \int_{S_R} u^2 \, dS\\
& \le & 2 C_3 R +   C_7\, \epsilon \int_{S_R} r^2 u^{{2n}\over
{n - 2}}\, dS 
+ C_8\, \epsilon \int_{B_o (R)} r u^{{2n}\over {n - 2}} \, dx + C_9
\epsilon^{-(n-2)/2} R \ \ \ \ \ \ \ \ \ \ 
\end{eqnarray*}
for large $R$, where we use Young's inequality and $C_7$, $C_8$ and
$C_9$ are positive constants. 
From  (1.9) we have 
$$
{d\over {dr}} \left( \int_{S_r} r^2 u^{{2n}\over {n - 2}} \, dS
\right) \le C_{10}  \int_{S_r} r u^{{2n}\over {n - 2}} \, dS \mfor r \ge
r_o\,,
$$
where $C_{10}$ is a positive constant. 
Hence 
\begin{eqnarray*}
(3.6) \ \ \ \ \ \ \ \ \ \ \ \ \ \ \ \ \  \int_{S_r} r^2 u^{{2n}\over {n - 2}}
\, dS & = & \int_0^r  
 \left( \int_{S_s} s^2 u^{{2n}\over {n - 2}} \, dS \right)' ds\\
& \le & C_{11} + \int_{r_o}^r 
 \left( \int_{S_s} s^2 u^{{2n}\over {n - 2}} \, dS \right)' ds\\
& \le & C_{11} + C_{10}\int_{r_o}^r \int_{S_s} s u^{{2n}\over {n - 2}} \, dS ds\\
& \le & C_{12} \int_{B_o (r)} r  u^{{2n}\over {n - 2}} dx \mfor r \ge
2 \,r_o\,,
\ \ \ \ \ \ \ \ \ \ \ \ \ \ \ \ \ \ \ \ 
\end{eqnarray*}
where $C_{11}$ and $C_{12}$ are  positive constants that depends on $r_o$
and $u$. If 
$$
\lim_{R \to \infty} \int_{B_o (R)}  u^{{2n}\over {n - 2}} \, dx <
\infty\,, \leqno (3.7)
$$
then by using theorem 3.1 in \cite{Leung.5}, which can be modified for $K$
satisfying (2.7) by making the constants there larger if necessary,
together with a result of Brezis and Kato \cite{Brezis-Kato}, we have 
$$
u (x) \le C_{13} |x|^{- (n - 2)}
$$
for large $|x|$ and for a positive constant $C_{13}$ (cf. the proof of
theorem 3.16 in \cite{Leung.5}). In particular, $u$
has slow decay. So we may assume that 
$$
\lim_{R \to \infty} \int_{B_o (R)}  u^{{2n}\over {n - 2}} \, dx = 
\infty\,. \leqno (3.8)
$$
As $K (x) \ge a^2$ for large $|x|$, it follows from (3.8) that
$$
\int_{B_o (R)} r K u^{{2n}\over {n - 2}} \, dx \ge {{a^2}\over 2} 
\int_{B_o (R)} r  u^{{2n}\over {n - 2}} \, dx \leqno (3.9)
$$
for large $R$. 
If we choose $\epsilon$ to be small, then (3.5), (3.6) and (3.9)
imply that 
$$
\int_{B_o (R)} r  u^{{2n}\over {n - 2}} \, dx \le C_{14} R \mfor R \ge r_o
$$
for large $R$ 
and
$$
\int_{S_R} R^2 u^{{2n}\over {n - 2}}\, dS \le C_{15} R 
$$
for large $R$, where $C_{14}$ and $C_{15}$ are positive constants.  
The last inequality can also be written as
$$
\int_{S^{n - 1}} u^{{2n}\over {n - 2}} (r, \theta)\, d\theta \le
C_{15} r^{-n} \leqno (3.10)
$$
for large $r$. Using lemma 4.39 and 4.40 in \cite{Leung.5}, we obtain 
spherical Harnack inequality (2.8) for spheres of large radius. Together
with (2.12) or (3.10) we conclude that $u$ has slow decay.\qedwh
{\bf Remark 3.11.} \ \ Condition (1.8) in theorem A is used only in obtaining a lower bound
on $P (u, r)$. Therefore theorem A remains valid if we replace (1.8)
by the assumption that $P (u, r) \ge - \delta^2$ for large $r$ and for
a constant $\delta$.

\vspace{0.5in}

{\bf \Large {\bf 4. \ \ Asymptotic behavior}}

\vspace{0.3in}

{\bf Theorem 4.1.} \ \ {\it  Assume that $\lim_{|x| \to \infty} K (x) =
K_\infty > 0$ and $|\btd
K| \le C_o$ in ${\R}^n$ for a positive number $C_o$. Let $u$ be
a positive smooth solution of equation (1.1) in ${\R}^n$ with slow decay. Assume
also that $K$ satisfies one of the conditions in lemma 2.4 so that $P
(u)$ exists. Then we have $P (u) \le 0\,.$ 
If $P (u) = 0$, then either $u$ has fast decay, or the conformal metric $u^{4/(n - 2)} g_o$ is complete 
and}
$$
\liminf_{|x| \to \infty} |x|^{{n - 2}\over 2} u (x) = 0
$$ 
{\it and there exists a sequence $\{s_l\} \subset \R$ of local maxima of $\bar v$, $s_l \to +\infty$ 
as $l \to +\infty\,,$ such that the 
sequence $\{v_l\}$ defined by}
$$
v_l (s, \,\theta) = v (s_l + s, \,\theta) \for s \in \R\,, \ \ \theta \in S^{n - 1} \ \ 
and \ \ l = 1, \ 2,..,
$$
{\it converges uniformly in $C^2$-norm on compact subsets of $\R \times S^{n - 1}$ to 
$C \,(\cosh s)^{(2 - n)/2}$
for a positive constant $C$ that depends on $K_\infty$. Here $v$ and $\bar v$ are  defined in
(1.11) and (1.12) respectively.}\\[0.1in] {\bf Proof.} \ \ Suppose that $P (u) > 0$. Lemma 2.13 
implies that there exist positive constants $c_1\,,$ $c_2$ and $s_o$ such that 
$$
c_1 \le v (s, \theta) \le c_2 \mfor s \ge s_o \ \ \ {\mbox{and}} \ \
\theta \in S^{n - 1}. \leqno (4.2)
$$
We have $P (u) = P (v) = \lim_{s \to + \infty} P (v, s)$. For $t > 0$,
let  
$$
v_t (s, \theta) = v (s + t, \theta) \mfor \theta \in S^{n - 1} \ \ \ \
{\mbox{and}} \ \ s \ge s_o - t\,.
$$
Using the equation
$$
{{\partial^2 v_t}\over {\partial s^2}} + \Delta_\theta \,v_t = F_t  \mfor s \ge
s_o - t\,, \leqno (4.3)
$$
where 
$$
F_t (s, \theta) = \left({{n
- 2}\over 2} \right)^2 v_t (s, \theta) - K (s + t, \theta) v_t^{{n + 2}\over {n - 2}}
(s, \theta)   \leqno (4.4)
$$
for $s \in \R$ and $\theta \in S^{n - 1}$, (4.2), the boundedness of $|\btd K|$, and elliptic estimates 
(cf. \cite{Gilbarg-Trudinger}), there exists a sequence of positive numbers
$t_i \to \infty$ as $i \to \infty$ and a $C^2$-function $v_\infty$
defined on $\R \times S^{n - 1}$ such that $v_{t_i}$ converges in
$C^2$-norm on compact subsets of $\R \times S^{n - 1}$ to
$v_\infty$ as $i \to \infty$. Furthermore, $v_\infty$ satisfies the equation
$$
{{\partial^2 v_\infty}\over {\partial s^2}} + \Delta_\theta\, v_\infty - \left({{n
- 2}\over 2} \right)^2 v_\infty + K_\infty\,  v_\infty^{{n + 2}\over {n
- 2}} = 0 \ \ \ \ {\mbox{in}} \  \ \R \times S^{n - 1}\,, \leqno (4.5)
$$
and $c_1 \le v_\infty \le c_2$ in $\R \times S^{n - 1}.$ 
By a result of Caffarelli, Gidas, Spruck in
\cite{Caffarelli-Gidas-Spruck}, 
$v_\infty$ is independent on $\theta$. 
It follows that the Pohozaev number $P
(v_\infty)$ is negative (cf. section 2.1 in \cite{K-M-P-S} ). Hence $P
(v_\infty, s) < 0$ for large $s$.  Fixed  a
large number $s_o$ so that  $P (v_\infty, s_o) < 0$. As $v_{t_i}$
converges in $C^2$-norm to $v_\infty$ in a compact neighborhood of $s_o \times S^{n
- 1}$, we have $P (v_{t_i}, s_o) < 0$ for $i$ large. But  $P (v_{t_i},
s_o)  = P (v, t_i + s_o) \to P(v) > 0$ as $\,i \to \infty\,$. Therefore we
have a contradiction. Hence $P (u) \le 0\,.$  Assume that $P (u) = 0$. 
It follows from the above argument that $v$ cannot be bounded away from zero in ${\R}^+ \times S^{n - 1}$, that is, 
$$
\liminf_{|x| \to \infty} |x|^{{n - 2}\over 2} u (x) = 0\,.
$$
Define 
$$
\bar v (s) = \int_{S^{n - 1}} v (s, \theta) \, d\theta \ \ \ \ \mfor s
\in {\R}\,. \leqno (4.6)
$$
Then $\bar v$ is bounded from above in ${\R}$. Either ${\bar v}' (s) \ge 0$ or ${\bar v}' (s) \le 0$ for all large
$s$, or ${\bar v}'$ keeps on changing sign near positive infinity. In the
first two cases we have $\lim_{s \to + \infty} \bar v (s) = \nu_o$
exists, as $\bar v (s)$ is bounded positive for large $s$. If $\nu_o >
0$, then using  spherical Harnack inequality (2.8), there
exists a positive constant $c_3$ such that 
$v (s, \theta) \ge c_3^2$ for large $s$ and for $\theta \in S^{n -
1}$. Hence we may apply the above argument to  show that $P (u) < 0$,
which is a contradiction. Assume that $\nu_o = 0$. 
Spherical Harnack inequality (2.8)  implies that 
$$
\lim_{s \to + \infty}\max_{\theta \in S^{n - 1}}  v (s, \theta) = 0\,.
$$
It seems standard to conclude that $u$ has fast decay. We show this in
the appendix. Therefore we need only to consider the case when  
${\bar v}'$ keeps on changing sign  near positive infinity. 
If this is the
case, then there exists a sequence of numbers $\{ s_i \}$ such that 
$\lim_{i \to \infty} s_i = + \infty$ and $\bar v$ achieves local
maximum at $s_i$ for $i = 1, 2,...$ From  (A.8) below, there exists a
positive constant $C_2$ such that 
$$
{\bar v}'' (s) - \left( {{n - 2}\over 2} \right)^2 {\bar v} (s) + C_2 \, 
{\bar v}^{{n + 2}\over {n - 2}} (s) \ge 0 \leqno (4.7)
$$
for large positive $s$. Using  (4.7) and the fact that ${\bar v}'' (s_i) \le 0$
we have 
$$
\bar v (s_i) \ge \left[ {1\over {C_2}} \left( {{n - 2}\over 2} \right)^2
\right]^{{n - 2}\over 4} \leqno (4.8)
$$
for large $i$. It follows from spherical Harnack inequality  (2.8)
that there exists a positive constant $c_4$ such that 
$$
v (s_i, \theta) \ge c_4^2 \mfor \theta \in S^{n - 1} \leqno (4.9)
$$
and for large $i$. Define
$$
v_i (s, \theta) = v (s + s_i\,, \,\theta) \mfor s \in \R\,, \ \ \theta \in S^{n - 1} 
\leqno (4.10)
$$
and $i = 1, 2,...$. As above,  
 a subsequence of $\{ v_i \}$  converges in
$C^2$-norm on compact subsets of $\R \times S^{n - 1}$ to a
non-negative $C^2$
function $v_\infty$ on $\R \times S^{n - 1}$
which satisfies the equation 
$$
{{\partial^2 v_\infty}\over {\partial s^2}} + \Delta_\theta v_\infty = v_\infty \left[ \left( {{n - 2}\over 2} \right)^2
- K_\infty \,
v_\infty^{4\over {n - 2}}  \right] \ \ \ \ {\mbox{in}} \ \ \R \times
S^{n - 1}. \leqno (4.11)
$$
Furthermore, from (4.9) we have
$$
v_\infty (0, \theta) \ge c_4^2 > 0 \mfor \theta \in S^{n - 1}.
$$   
Let 
$$
u_\infty (x) = u_\infty (r, \theta) = r^{- (n - 2)/2} \, v_\infty (s, \theta) \mfor |x| >
0\,,
$$
where $x = (r, \theta)$ and $r = e^s$. Then $u_\infty$ satisfies the
equation 
$$
\Delta u_\infty = - K_\infty\, u_\infty^{{n + 2}\over {n - 2}} \ \ \ \
{\mbox{in}} \ \ {\R}^n \setminus \{ 0 \}\,. \leqno (4.12)
$$
Furthermore $u_\infty (x) \ge 0$ for $|x| > 0$ and $u_\infty (x) \ge c_4^2 >
0$ for $|x| = 1$. It follows from equation (4.12) and the maximum
principle that $u_\infty (x) > 0$ for $|x| > 0$. Hence $v_\infty$ is
positive in $\R \times S^{n - 1}$. Because $v_\infty$ may not be bounded away from zero, we 
can only concluded that $P (v_\infty) \le 0$. The above argument shows that if 
$P (v_\infty) < 0$, then $P (u) < 0$. Hence $P (v_\infty) = 0\,.$  
By a result of Caffarelli, Gidas and Spruck (\cite{Caffarelli-Gidas-Spruck},  
see also \cite{K-M-P-S}),  $v_\infty$ is radial, and together with $v_\infty' (0) = 0$ we have 
$v_\infty (s) = C (\cosh s)^{(2-n)/2}$ for $s \in \R$ and for a positive number $C$ depending on $K_\infty$. 
It remains to show that the conformal metric corresponding to 
$u$ is complete in this case. From (4.7) and the fact that $v$ is bounded from above, 
there exists a positive constant $c_5$ such that 
$$
{\bar v}'' (s) \ge - c_5^2 \mfor s \in \R\,. \leqno (4.13)
$$
At $s_i$ we have ${\bar v}' (s_i) = 0$, therefore ${\bar v}' (s)$ is not too negative for $s$ close to $s_i$ 
for $i =1, 2,...$ By using (4.8) and (4.13), there exists a positive number $\epsilon$ independent  
on $i$ such that 
$$
{\bar v} (s) \ge {1\over 2}  \left[ {1\over {C_2}} \left( {{n - 2}\over 2} \right)^2
\right]^{{n - 2}\over 4} \mfor s \in [s_i, \,s_i + \epsilon]\,, \ \ i = 1, 2,... \leqno (4.14)
$$
Using spherical Harnack inequality (2.8) there exists a positive constant $c_6$ 
such that 
$$
v (s, \theta) \ge c_6^2 \mfor s \in [s_i, \,s_i + \epsilon]\,, \ \ \theta \in S^{n - 1} \leqno (4.15)
$$
and $i = 1, 2,...$ Without loss of generality, we may assume that $s_1 > 0$ and $s_{i + 1} > s_i$ for 
$i = 1, 2,...$ 
For any fixed $\theta \in S^{n - 1}$, the length of the curve $r \mapsto (r, \theta)$, $r \in (1, \infty)$,  
in the conformal metric $u^{4/(n - 2)} g_o$ is given by 
$$
\int_1^\infty u^{2\over {n - 2}} (r, \theta) \, dr = \int_0^\infty v^{2\over {n - 2}} (s, \theta) \, ds 
\ge \sum_{i = 1}^\infty \int_{s_i}^{s_i + \epsilon} v^{2\over {n - 2}} (s, \theta) \, ds  = \infty\,.
$$
Hence the conformal metric is complete. \qedwh
{\bf Corollary 4.16.} \ \ {\it Assume that $K$ satisfies (2.7). Let $u$ be
a positive smooth solution of equation (1.1) in ${\R}^n$ with slow decay. Then either $u$ has fast 
decay or the conformal metric $g = u^{4/(n - 2)} g_o$ is complete. Furthermore,  
the conformal metric $g$ is complete  if and only if the total volume of $({\R}^n, g)$ is infinite.}\\[0.1in]
{\bf Proof.} \ \ Consider $\bar v$ as defined in (4.6) and $\lim_{s \to \infty} {\bar v} (s)$. If
the limit exists and is equal to zero, then $u$ has fast decay (see appendix A). If the limit is non-zero, then 
$u (x) \ge c |x|^{(2 - n)/2}$ for large $|x|$ and for a positive constant $c$, and hence the  
conformal metric $g$ is complete. In case the limit does not exists, then ${\bar v}'$ 
keeps on changing signs near positive infinity. The argument in the last part of the proof of theorem 
4.1 shows that the conformal metric is also complete in this case. Finally, if the total volume of 
$({\R}^n, g)$ is infinite, then $u$ does not have fast decay and hence by above $g$ is complete. If
$g$ is complete, then either $\lim_{s \to + \infty} \bar v (s)$ exists and is positive, or (4.15) holds. In either cases we have
$$
\int_{{\R}^n \setminus B_o (1)} u^{{2n}\over {n - 2}} (x) \, dx = \int_0^\infty \int_{S^{n - 1}} 
v^{{2n}\over {n - 2}} \, d\theta \, ds = \infty\,.
$$
That is, the total volume of $({\R}^n, g)$ is infinite.\qedwh
{\bf Theorem 4.17.} \ \ {\it Assume that $K$ satisfies the conditions in theorem 4.1 and $x \cdot
\btd K (x) \le 0$ for large $|x|$. 
 Let $u$ be
a positive smooth solution of equation (1.1) in ${\R}^n$ with slow decay. If $P (u) = 0$, then $u$ has fast
decay.}\\[0.1in]
{\bf Proof.} \ \ Let $\bar v$ be defined as in (4.6). Suppose that $P (u) = 0$
and $u$ does not have fast decay. It follows from the argument in the proof of theorem 4.1 that there
exists a sequence $\{s_j\} \subset \R$ such that $s_j \to + \infty$ as $j \to \infty$ and each $s_j$
is a local {\it minimum} for $\bar v$ and $\lim_{j \to \infty} \bar v (s_j) = 0$. Let 
$$
w_j \,(s, \theta) = {{v (s + s_j, \,\theta)}\over {\bar v (s_j)}} \mfor s \in R\,, \ \ \theta \in S^{n
- 1} \ \ {\mbox{and}} \ \ j = 1, 2,... \leqno (4.18)
$$
Given $S > 0$, we claim that there exist positive numbers $C_s$ and $k_s$ such that 
$$
w_j \,(s, \theta) \le C_s \mfor s \in [- S, S]\,, \ \ \theta \in S^{n
- 1} \ \ {\mbox{and}} \ \ j \ge k_s\,. \leqno (4.19)
$$
From (A.7) and (A.8) in the appendix we have  
$$
\left( {{n - 2}\over 2} \right)^2 {\bar v} (s) - C_1 \,
{\bar v}^{{n + 2}\over {n - 2}} (s) 
\ge \, {\bar v}'' (s) \, \ge \left( {{n - 2}\over 2} \right)^2 {\bar v} (s) - C_2 \,
{\bar v}^{{n + 2}\over {n - 2}} (s)   
$$
for large $s$ and for some positive constants $C_1$ and $C_2$. It follows that there exists a
positive constant $\varepsilon_o$ such that 
$$
0 \le {\bar v}'' (s) \le {1\over 2} \left( {{n - 2}\over 2} \right)^2 {\bar v} (s)
\ \ \ \ {\mbox{for}} \ \ \bar v (s) \le \varepsilon_o \ \ {\mbox{and \ \ large}} \ \ s\,. \leqno
(4.20)
$$
Choose the number $k_s$ such that  
$$
\bar v (s_j) \, e^{ {1\over 4} \left( {{n - 2}\over 2} \right)^2 S^2}  < \varepsilon_o \mfor j \ge
k_s\,.
$$
We may also assume that for $j \ge k_s$, $s_j$ is large. 
For $j \ge k_s$ and $s$ close to $s_j\,,$ $\,{\bar v}'' (s)$ satisfies (4.20). Hence we obtain
$$
{\bar v}' (s) \le {1\over 2} \left( {{n - 2}\over 2} \right)^2 \int_{s_j}^s \bar v (t) \, dt \le 
{1\over 2} \left( {{n - 2}\over 2} \right)^2 \bar v (s) (s - s_j) \leqno (4.21)
$$
for $s > s_j$ close to $s_j\,,$ as $v' (s) \ge 0$ for $s \ge s_j$  close to $s_j$. Therefore
we have
$$
{{\bar v (s)}\over { \bar v (s_j)}} \, \le \,  \, e^{ {1\over 4} \left( {{n - 2}\over 2} \right)^2 (s
- s_j)^2}
\
\
\
\
\mfor s \in [s_j, \,s_j + S]\,.
$$
Likewise, we obtain a similar inequality on $[s_j - S, \,s_j]$.  By using spherical Harnack inequality
(2.8), we obtain the desired constant $C_s$ and the bound in (4.19). It follows from (4.19) that
there exists a subsequence of
$\{w_j\}$ which converges in $C^2$-norm on compact subsets of $\R \times S^{n - 1}$ to a
solution $w$ of the equation 
$$
{{\partial^2 w}\over {\partial s^2}} + \Delta_\theta \,w - \left( {{n - 2}\over 2} \right)^2 w = 0
\ \ \ \ {\mbox{in}} \ \ \R \times S^{n - 1}.
$$
The associated function $h$ related to $w$ as in (1.11) is, by the maximum principle and the above
equation, a positive harmonic function on ${\R}^n \setminus \{0\}.$ Therefore $h (x) = a |x|^{2 - n}
+ b$ and 
$$
w (s, \theta) = a \,e^{- {{n - 2}\over 2} s} + b\,e^{{{n - 2}\over 2} s} \mfor s \in \R \ \
{\mbox{and}}
\ \ \theta \in S^{n - 1}  \leqno (4.22) 
$$
for some positive constants $a$ and $b$. As $\bar w$ has a critical point at $s = 0$, 
we have $a = b
> 0$. As in
\cite{K-M-P-S} we obtain
\begin{eqnarray*}
(4.23) & \ & \lim_{j \to \infty} \int_{S^{n - 1}} \left[ {1\over 2} 
\left( {{\partial w_j}\over {\partial s}} \right)^2  (0, \theta) 
-  {1\over 2} | \btd_\theta w_j |^2 (0, \theta)
- {1\over 2} \left( {{n - 2}\over 2} \right)^2  w_j^2 (0, \theta)
\right. \ \ \ \ \ \ \ \ \ \ \ \ \ \ \  \\
&\ & \ \ \ \ \ \ \ \ \ \ \ \ \ \ \left. \, + \, {\bar v}^{4\over {n - 2}} (s_j) \,{{n - 2}\over
{2n}}  K  \,w_j^{{2n}\over {n - 2}} (0, \theta) \right] \, d\theta \ \ \ \ \ \ \ \\
& = & \int_{S^{n - 1}} \left[ {1\over 2} 
\left( {{\partial w}\over {\partial s}}   \right)^2 (0, \theta)  
- {1\over 2} \left( {{n - 2}\over 2} \right)^2  w^2 (0, \theta) \right] \, d\theta\\
& = & -\omega_n \left({{n - 2}\over 2} \right)^2 ab < 0\,. 
\end{eqnarray*}
Using the assumption that $({\partial K}/\partial s) (s, \theta) \le 0$ for large $s$ and
$\theta
\in S^{n - 1},$  Pohozaev identity (2.25) and 
$\lim_{s \to + \infty} P (v, s) = 0$, we see that $P (v, s) \ge 0$ for large $s$. On
the other hand, by (4.23) we have
\begin{eqnarray*}
 {{P (v, s_j)}\over {{\bar v}^2 (s_j)}} 
& = &\int_{S^{n - 1}} \left[ {1\over 2} 
\left( {{\partial w_j}\over {\partial s}}  \right)^2  (0, \theta)
-  {1\over 2} | \btd_\theta w_j |^2 (0, \theta)
- {1\over 2} \left( {{n - 2}\over 2} \right)^2  w_j^2 (0, \theta)
\right. \ \ \ \ \ \ \ \ \\
& \ &  \ \ \ \ \ \ \ \ \ \ \ \  \left. \, + \, {\bar v}^{4\over {n - 2}} (s_j) {{n - 2}\over {2n}}  K 
w_j^{{2n}\over {n - 2}} (0, \theta)\right] \, d\theta < 0  \ \ \ \ \ \ \  \ \ \ \ \ \ \ 
\end{eqnarray*}
for $j$ large. Hence we have a contradiction. Therefore if $P (u) = 0$, then $u$ has fast
decay.\qedwh
{\bf Remark 4.24.} \ \ One can also prove theorem 4.17 by modifying the argument of Chen and Lin 
in section 3 of
\cite{Chen-Lin.1}.\\[0.2in]
{\bf Theorem 4.25. } \ \ {\it Assume that $K$ satisfies (2.7). Let $u$ be
a positive smooth solution of equation (1.1) in ${\R}^n$. Assume also that }
$$
\lim_{r \to \infty}  \int_{S_r} r \left[ {{\partial u}\over {\partial r}} + {{n
- 2}\over 2} {u\over r} \right]^2 \, dS = 0\,. \leqno (4.26)
$$
{\it If $P (u)$ exists, then $P (u) \le 0$. Furthermore, $P (u) = 0$
implies that}
$$
\,\liminf_{|x| \to \infty} \,|x|^{(n - 2)/2} \,u (x) = 0\,. \leqno (4.27)
$$

\vspace*{0.05in}

{\bf Proof.} \ \ Let 
$$
\omega (r) = r^{n - 2} \int_{S^{n - 1}} u^2 (r, \theta)\, d\theta \mfor
r > 0\,. \leqno (4.28)
$$
We have 
$$
\omega' (r) = 2\, r^{n - 2} \int_{S^{n - 1}} u \left[ {{\partial u}\over
{\partial r}} + {{n - 2}\over 2} {u \over r} \right] \,d \theta \mfor r > 0\,.
\leqno (4.29)
$$
Using equation (1.1) we obtain
$$
{{\omega'' (r)}\over 4}  = {{n - 3}\over 2} r^{n - 3} \int_{S^{n - 1}} u \left[
{{\partial u}\over {\partial r}} + {{n - 2}\over 2} {u \over r}
\right] \, d\theta\, +\, {{r^{n - 2}}\over 2} \left[ \int_{S^{n - 1}}
|\btd u|^2 \, d\theta - \int_{S^{n - 1}} K u^{{2n}\over {n - 2}} \, d
\theta \right] \leqno (4.30)
$$
for $r > 0\,.$ From Pohozaev identity (2.2) we also have
\begin{eqnarray*}
(4.31) \ \ \ \ 
P (u, r) & = & \int_{S_r} r \left[ {{\partial u}\over {\partial r}} + {{n
- 2}\over 2} {u\over r} \right]^2 \, dS - {{n - 2}\over 2} 
 \int_{S_r} u \left[ {{\partial u}\over {\partial r}} + {{n
- 2}\over 2} {u\over r} \right] \, dS\ \ \ \ \ \ \ \ \ \ \ \ \ \ \ \ 
 \ \ \ \ \ \ \ \ \ \\
& \ & \ \ \ \ \ \ \  - {1\over 2} \int_{S_r} r | \btd
u |^2 \,dS + {{n - 2}\over {2n}} \int_{S_r} r K u^{{2n}\over {n - 2}}
\, dS
\end{eqnarray*}
for $r > 0\,.$ 
Under condition (4.26) and $P (u) \ge 0$, we claim that $\omega' (r)$
cannot be non-negative for all large $r$. In order to obtain a
contradiction, suppose that $\omega' (r) \ge 0$ for large $r$.  
From (4.29) we have 
$$
\omega' (r) \le  2 \, \left( \int_{S_r} {{u^2}\over r} \, dS
\right)^{1\over 2} 
\left( \int_{S_r}  {1\over r}  \left[ {{\partial u}\over {\partial r}} + {{n
- 2}\over 2} {u\over r} \right]^2 \,dS \right)^{1\over 2}
 \le {{C_1\, \omega^{1\over 2} (r)}\over r} \leqno (4.32)
$$
for large $r$, where $C_1$ is a positive constant. If there exists a
positive constant $c'$ such that $\omega' (r) \ge c'$ for large $r$,
then we have 
$$
\omega (r) \ge {{c' \, r}\over 2} \leqno (4.33)
$$
for large $r$. It follows from (4.32) that 
$$
 \omega' (r) \le {{C_2\, \omega (r)}\over {r^{3\over 2}}} \leqno (4.34)
$$
for large $r$ and a positive constant $C_2$. Integrating both sides of (4.34) we conclude that
$\omega$ is bounded from above in ${\R}^+$, which contradicts
(4.33). Therefore we may assume that 
$$
\liminf_{r \to \infty}\, \omega' (r)
= 0\,. \leqno (4.35)
$$ 
Let $c$ to be a positive constant to be fixed later.   
Suppose that the inequality
$$
\omega'' (r) \le - c \,{{\omega (r)}\over {r^2}} \leqno (4.36) 
$$
does not hold for all large $r$. Then there exists a sequence of positive numbers
$r_j \to \infty$ as $j \to \infty$ such
that 
$$
\omega'' (r_j) > - c \,{{\omega (r_j)}\over {r_j^2}} \leqno (4.37)
$$
for $j = 1, 2,...\,$ From (4.30) we have
$$
{1\over 2} \int_{S_{r_j}} r_j  |\btd u|^2 \, dS + {{n - 3}\over 2} 
\int_{S_{r_j}} u  \left[ {{\partial u}\over {\partial r}} + {{n
- 2}\over 2} {u\over {r_j}} \right] \, dS + {{c\, \omega (r_j)}\over 4} > 
{1\over 2} \int_{S_{r_j}} r_j K u^{{2n}\over {n - 2}} \, dS \leqno
(4.38)
$$
for $j = 1, 2,...$ On the other hand, from (4.31) we have
\begin{eqnarray*}
(4.39) \ \ \ \ & \ & 
{1\over 2} \int_{S_{r_j}} r_j K  u^{{2n}\over {n - 2}} \, dS  =  P (u, r_j)
+ {1\over 2} \int_{S_{r_j}} r_j  |\btd u|^2 \, dS\\
& \ & \ \ \ \  + \,{{n - 3}\over 2} 
\int_{S_{r_j}} u  \left[ {{\partial u}\over {\partial r}} 
+ {{n - 2}\over 2} {u\over {r_j}} \right] \, dS  +\, {1\over 2}\int_{S_{r_j}} u  \left[ {{\partial u}\over {\partial r}}
+ \,{{n - 2}\over 2} {u\over r_j} \right] \, dS\\
& \ & \ \ \ \  
 - \, 
\int_{S_{r_j}} r_j \left[ {{\partial u}\over {\partial r}} 
+ {{n - 2}\over 2} {u\over r_j} \right]^2 \, dS  +\, {1\over n} 
\int_{S_{r_j}} r_j K u^{{2n}\over {n - 2}} \, dS \ \ \ \ \ \ \ \ \ \ \
\ \ 
\end{eqnarray*}
for $j = 1, 2,...$ Using the H\"older inequality and the Young inequality we obtain
$$
\omega (r_j) = \int_{S_{r_j}} {{u^2}\over {r_j}}\, dS \le
\omega_n^{2/n} \left( \int_{S_{r_j}} r_j \, u^{{2n}\over {n - 2}} \,
dS \right)^{{n - 2}\over n} \leqno (4.40)
$$
and 
$$
\omega (r_j) \le \,{{n - 2}\over n} \int_{S_{r_j}} r_j \,
u^{{2n}\over {n - 2}} \,dS +  {2 \over n} \omega_n \leqno (4.41)
$$
for $j = 1, 2,...\,$ As $\omega' (r) \ge 0$ for large $r$ and $K (x)
\ge a^2$ for large $|x|$, it follows from (4.40) that  there exists
a positive number $c_o$ such that 
$$
 \int_{S_{r_j}} r_j K \, u^{{2n}\over {n - 2}} \,
dS \ge c_o^2 \leqno (4.42)
$$
for $j = 1, 2,....$ Using (4.26), (4.29), (4.41) and (4.42) we have
\begin{eqnarray*}
(4.43) \ \ {1\over n} 
\int_{S_{r_j}} r_j K u^{{2n}\over {n - 2}} + P (r, r_j) & + & 
{1\over 2}\int_{S_{r_j}} u  \left[ {{\partial u}\over {\partial r}} 
+ {{n - 2}\over 2} {u\over r_j} \right] \, dS\\
& \ & \ \ \ - 
\int_{S_{r_j}} r_j \left[ {{\partial u}\over {\partial r}} 
+ {{n - 2}\over 2} {u\over r_j} \right]^2 \, dS > {{c \omega
(r_j)}\over 4} \ \ \ \ \ \ \ \ 
\end{eqnarray*}
for large $j$, if we choose $c$ to be small enough. (4.38), (4.39)
together with (4.43) give a contradiction. Hence (4.36) holds for a small
$c$. Integrating both sides of (4.36) and using the assumption that
$\omega (r)$ is non-decreasing  for large $r$, we obtain
$$
\omega' (R) - \omega' (r) \le - c \,\int_r^R {{\omega (s)}\over {s^2}} \,
ds \le - c\, \omega (r) \int_r^R {{ds}\over {s^2}} = c\, \omega (r)
\left[ {1\over R} - {1\over r} \right] \leqno (4.44)
$$
for large $R$ and $r$ with  $R > r$. Letting $R$ approaches infinity by a sequence
which satisfies (4.35) we have
$$
\omega' (r) \ge {{c\, \omega (r)}\over r} \mfor r \ge r_o\,. \leqno (4.45)
$$
Integrating both sides of (4.45), there exists a positive constant $C_3$
such that 
$$
\omega (r) \ge C_3 \, r^c \mfor r \ge r_o\,. \leqno (4.46)
$$
Together with (4.32) we have
$$
\omega' (r) \le {{C_4 \, \omega (r)}\over {r^{1 + {c\over 2}} }}\leqno
(4.47)
$$
for large $r$, where $C_4$ is a positive constant. Integrating both
sides of (4.47) we conclude again that $\omega$ is bounded from above in
${\R}^+$, which contradicts (4.46). Hence we may dismiss the case that $\omega' (r)$ is
non-negative for all large $r$.\bk
Next, assume that $\omega$ achieves a local minimum at $r_o$. Therefore
we have $\omega' (r_o) = 0$ and $\omega'' (r_o) \ge 0$. It follows
from (4.29), (4.30) and (4.31) that 
\begin{eqnarray*}
& \ & {1\over 2} \int_{S_{r_o}} r_o  |\btd u|^2 \, dS  \ge 
{1\over 2} \int_{S_{r_o}} r_o \,K u^{{2n}\over {n - 2}} \, dS 
\ \ \ \ {\mbox{and}}\\
{1\over 2} \int_{S_{r_o}} r_o \,K u^{{2n}\over {n - 2}} \, dS & =&  P (u, r_o)
+ {1\over 2} \int_{S_{r_o}} r_o  |\btd u|^2 \, dS
  - \int_{S_{r_o}} r_o \left[ {{\partial u}\over {\partial r}} 
+ {{n - 2}\over 2} {u\over r_o} \right]^2 \, dS\\
& \ & \ \ \ \   + {1\over n} 
\int_{S_{r_o}} r_o \,K u^{{2n}\over {n - 2}} \, dS\,. 
\end{eqnarray*} 
Therefore if $P (u) > 0\,,$ or $P (u) = 0$ and $\omega (r) \ge c_1^2 >
0$ for large $r$ and for a positive constant $c_1$, then $\omega$ has
no local minimum on $(r_1\,, \infty)$ for large $r_1$. Hence $\omega'
(r) \le 0$ for large $r$ and $\omega (r) \to C_o \ge 0$ as $r \to
\infty$. Using a result in \cite{Ding-Ni}, there exist a sequence $\{r_i\}$ of
positive numbers such that $\lim_{i \to \infty} r_i = \infty\,,$
$$
r_i\, \omega' (r_i) \to 0 \ \ \ \ {\mbox{and}} \ \ \ \ r_i\, [r_i \omega''
(r_i) + \omega' (r_i) ] \to 0 \ \ \ \ {\mbox{as}} \ \ i \to \infty\,.
\leqno (4.48)
$$
The first limit in (4.48) shows that 
$$
\int_{S_{r_i}} u \left[ {{\partial u}\over {\partial r}} + {{n
- 2}\over 2} {u\over r_i} \right] \, dS \to 0 \ \ \ \ 
{\mbox{as}} \ \ i \to \infty\,.
\leqno (4.49)
$$
The second limit in (4.48) together with (4.49) imply that 
$$
 \int_{S_{r_i}} r_i  |\btd u|^2 \, dS - \int_{S_{r_i}} r_i K
u^{{2n}\over {n - 2}} \, dS \to 0\ \ \ \ 
{\mbox{as}} \ \ i \to \infty\,.
\leqno (4.50)
$$
It follows from (4.26), (4.31), (4.49) and (4.50) that 
$$
P (u, r_i) \to - {1\over n} 
\int_{S_{r_i}} r_i K u^{{2n}\over {n - 2}} \, dS \le 0\ \ \ \ 
{\mbox{as}} \ \ i \to \infty\,.
\leqno (4.51)
$$
Hence $P (u) \le 0$, which is a contradiction. From (4.40) and (4.51), 
if $\omega (r) \ge c_1^2$ for large $r$, then $P (u) <
0$. Therefore if $P (u) = 0$, then $\liminf_{r \to \infty} \omega (r) =
0$, which in turns implies that $\liminf_{|x| \to \infty} |x|^{(n -
2)/2} u (x) =
0\,.$ \qed
Let $u_o = C\, r^{- (n - 2)/2}$ for $r \ge 1$, where $C$ is a positive
constant.  $u_o$ satisfies equation (1.1) in ${\R}^n
\setminus B_o (1)$ with $K = [(n - 2)/2]^2 \,C^{-4/(n - 2)}$. Furthermore, 
$$
{{\partial u_o}\over {\partial r}} + {{n - 2}\over 2} {{u_o}\over r} =
0 \ \ \ \ {\mbox{in}} \ \ {\R}^n
\setminus B_o (1)\,.
$$
Therefore condition (4.26) can also be viewed as $u$ is asymptotically
close to $u_o$. The results obtained here suggest that equation (1.1)
favors $x \cdot \btd K (x)$ to be negative rather than
positive. Indeed, Bianchi shows that if $K$ is radially 
symmetric and $K' \ge
0\,, \ \not\equiv 0$ in ${\R}^+$,  then equation (1.1) has no positive solutions in
${\R}^n$ \cite{Bianchi.New}.

\vspace*{0.5in}

{\bf \Large {\bf Appendix}}

\vspace*{0.3in}

{\bf Theorem A.1.} \ \ {\it Assume that $a^2 \le K (x) \le b^2$ for large $|x|$ and for positive 
constants $a$ and $b$. Let $u$ be a positive smooth solution of
equation (1.1) in
${\R}^n$. If $\,\lim_{|x| \to \infty} |x|^{(n - 2)/2} u (x) = 0\,,$ then
there exist positive constants $c_1$ and $c_2$ such that $c_1\, |x|^{2 - n} \le u (x) \le c_2\, 
|x|^{2 - n}$ for large $|x|$. In addition, if $u$ is radial, then  
$\lim_{r
\to
\infty} r^{n - 2} u (r)$ exists and is positive.}\\[0.1in] {\bf Proof.} \ \ Let 
$$
\bar u (r) = \int_{S^{n - 1}} u (r, \theta)\, d\theta \mfor r >
0\,. \leqno (A.2)
$$
We have 
$$
\int_{S^{n - 1}}K (r, \theta) \, u^{{n + 2}\over {n - 2}} (r, \theta)\, d\theta
\ge a^2\omega_n^{-{4\over {n+2}}}\, {\bar u}^{{n + 2}\over {n -
2}} (r) \leqno (A.3)
$$
for large $r$. It follows from equation (1.1) that
$$
{\bar u}'' (r) + {{n - 1}\over r} {\bar u}' (r) + C_1 \,{\bar u}^{{n +
2}\over {n - 2}} \le 0 \leqno (A.4)
$$
for large $r$, where $C_1 = a^2\omega_n^{-4/(n+2)}$ is a positive
constant.  Using spherical Harnack inequality (2.8) we have
\begin{eqnarray*}
\int_{S^{n - 1}}K (r, \theta)  u^{{n + 2}\over {n - 2}} (r, \theta)\,
d\theta & \le & b^2 \omega_n \left[ \sup_{S_r} u \right]^{{n + 2}\over
{n - 2}} = b^2\omega_n^{-{4\over {n-2}}}
 \left[\omega_n \sup_{S_r} u \right]^{{n + 2}\over
{n - 2}}\\
& \le & b^2\omega_n^{-{4\over {n-2}}}
 \left[C_h \int_{S^{n - 1}} u (r, \theta)\, d\theta \right]^{{n + 2}\over
{n - 2}}
\end{eqnarray*}
for large $r$. Therefore we obtain 
$$
{\bar u}'' (r) + {{n - 1}\over r} {\bar u}' (r) + C_2 \,{\bar u}^{{n +
2}\over {n - 2}} \ge 0 \leqno (A.5)
$$
for large $r$, where $C_2 = b^2\, C_h^{{n + 2}\over {n - 2}}
\omega_n^{-{4\over {n-2}}}$. Let 
$$
{\bar v} (s) = r^{{n - 2}\over 2} \bar u (r)\,, \ \ \ \ {\mbox{where}}
\ \ \ \ r = e^s\,, \ \ r >>1\,. \leqno (A.6)
$$
From (A.4) and (A.5) we have
\begin{eqnarray*}
(A.7)  & \ & \ \ \ \ \ \ \ \ \ \ \ \ \ \ 
{\bar v}'' (s) - \left( {{n - 2}\over 2} \right)^2 {\bar v} (s) + C_1 \,
{\bar v}^{{n + 2}\over {n - 2}} (s) \le 0 \ \ \ \ {\mbox{and}} \ \ \ \ \ \ \ \ \ \
\ \ \ \ \ \ \ \ \ \ \ \ \ \ \ 
\\
(A.8) & \ & \ \ \ \ \ \ \ \ \ \ \ \ \  
{\bar v}'' (s) - \left( {{n - 2}\over 2} \right)^2 {\bar v} (s) + C_2 \, 
{\bar v}^{{n + 2}\over {n - 2}} (s)  \ge  0
\end{eqnarray*}
for large positive $s$. As ${\bar v} (s) \to 0$ when $s \to + \infty$, (A.8)
implies that ${\bar v}'' (s) > 0$ for large $s$. Hence ${\bar v}' (s) < 0$
for large $s$. That is,
$$
- r\, {\bar u}' (r) > {{n - 2}\over 2} \,\bar u (r) \leqno (A.9)
$$
for large $r$. Multiplying both sides of (A.9) by $C_2\,r^{-1} [2/(n - 2)]\, [r^{(n -
2)/2} {\bar u} (r) ]^{4/(n - 2)}$ and using the fact that 
$r^{(n - 2)/2}\, {\bar u} (r) \to 0$ as $r \to \infty$, for any
$\varepsilon > 0$, we have
$$
- \varepsilon \, {\bar u}' (r) > C_2\, r \, {\bar u}^{{n + 2}\over {n - 2}} (r)
\leqno (A.10)
$$
for large $r$. As ${\bar u}' (r) \le 0$ for large $r$, from (A.10) we
obtain 
$$
\varepsilon \, {\bar u} (r) > C_2 \int_r^\infty t \,{\bar u}^{{n + 2}\over {n
- 2}} (t) \, dt \leqno (A.11)
$$
for large $r$. From (A.5) we have
$$
[r {\bar u}' (r)]' + (n - 2) {\bar u}' (r) + C_2 \, r \, {\bar u}^{{n
+ 2}\over {n - 2}} (r)  \ge 0  \leqno (A.12)
$$
for large $r$. Integrating both sides of (A.12) from $r$ to $R$, letting
$R \to \infty$ and using the gradient estimate in  lemma 2.9, we obtain
$$
r {\bar u}' (r) + (n - 2)\, {\bar u} (r) \le C_2 \int_r^\infty t \,{\bar u}^{{n + 2}\over {n
- 2}} (t) \, dt < \varepsilon \, {\bar u} (r) \leqno (A.13)
$$
for large $r$. It follows from (A.13) that for any positive number $m
< (n - 2)$, there exists a positive constant $C_3$ such that 
$$
{\bar u} (r) \le C_3 \,r^{-m} \leqno (A.14)
$$
for large $r$. Similarly we have
$$
r {\bar u}' (r) + (n - 2)\, {\bar u} (r) \ge C_1 \int_r^\infty t \,{\bar u}^{{n + 2}\over {n
- 2}} (t) \, dt \leqno (A.15)
$$
for large $r$. It follows from (A.4) that
$$
 r^{n - 1}  {\bar u}' (r)  \le r_o^{n - 1} {\bar u}' (r_o) 
- C_1 \, \int_{r_o}^r t^{n - 1} {\bar u}^{{n +
2}\over {n - 2}} (t) \, dt \leqno (A.16)
$$
for $r_o$ and $r$ large, with $r > r_o$. We note that the integral in
(A.16) makes sense because of (A.14). From (A.15) we have
$$
 r^{n - 1}  {\bar u}' (r) \ge  C_1 \,r^{n - 2} 
\int_r^\infty t \,{\bar u}^{{n + 2}\over {n
- 2}} (t) \, dt - (n - 2) r^{n - 2} {\bar u} (r) \leqno (A.17)
$$
for large $r$. Fix $r_o$ to be large enough so that by (A.9) $u' (r_o)
\le 0$. Substituting (A.17) into (A.16) we obtain
$$
r^{n - 2} {\bar u} (r) \ge {1\over {n - 2}} \left[ 
 C_1 \, \int_{r_o}^r t^{n - 1} {\bar u}^{{n +
2}\over {n - 2}} (t) \, dt +  C_1 r^{n - 2} 
\int_r^\infty t \,{\bar u}^{{n + 2}\over {n
- 2}} (t) \, dt -r_o^{n - 1} {\bar u}' (r_o) \right] \ge C_4^2
$$
for large $r$ and for some positive constant $C_4$.  Likewise,
using (A.5) and (A.13) we have
$$
r^{n - 2} {\bar u} (r) \le {1\over {n - 2}} \left[ 
 C_2 \, \int_{r_o}^r t^{n - 1} {\bar u}^{{n +
2}\over {n - 2}} (t) \, dt +  C_2 r^{n - 2} 
\int_r^\infty t \,{\bar u}^{{n + 2}\over {n
- 2}} (t) \, dt -r_o^{n - 1} {\bar u}' (r_o) \right] \le C_5^2
$$
for large $r$ and for some positive constant $C_5$. Using 
spherical Harnack inequality (2.8), we obtain the desired bounds. We
observe that in case
$K$ and $u$ are radially symmetric, the above argument shows that
$$
r^{n - 2} u (r) = {1\over {n - 2}} \left[ 
  \int_{r_o}^r K (t)\,t^{n - 1} \,
 u^{{n +
2}\over {n - 2}} (t) \, dt +   r^{n - 2} 
\int_r^\infty K (t) \,
 t \, u^{{n + 2}\over {n
- 2}} (t) \, dt -r_o^{n - 1} u' (r_o) \right]  \leqno (A.18)
$$
for large $r$. Hence
$\lim_{r \to \infty} r^{n - 2}\, u (r) = c_o > 0$, where 
$$
c_o =  {1\over {n - 2}} \left[ 
  \int_{r_o}^\infty K (t)\,t^{n - 1} \,
 u^{{n +
2}\over {n - 2}} (t) \, dt -r_o^{n - 1} u' (r_o) \right] 
$$
by using (A.18).\qed
It follows from lemma A.1 and (2.12) that the only possible indices
$\kappa$ such that 
$$
c_1 \le |x|^\kappa u (x) \le c_2 \ \ \ \ {\mbox{for \ \ large}} \ \ |x|
$$
and for some positive numbers $c_1$ and $c_2$ are $\kappa = (n - 2)/2$ or
$\kappa = n - 2$.\\[0.2in]

\vspace*{0.5in}

%\pagebreak

\end{document}